\documentclass[final,leqno]{siamltex}

\usepackage{booktabs,caption,fixltx2e}
\usepackage[flushleft]{threeparttable}
\usepackage{amsmath, amssymb}
\usepackage{latexsym, color}
\usepackage{epsfig}
\usepackage{bm}
\usepackage{tikz}

\numberwithin{equation}{section}

\newtheorem{thm}{Theorem}[section]

\newcommand{\eqnref}[1]{(\ref {#1})}

\newcommand{\beq}{\begin{equation}}
\newcommand{\eeq}{\end{equation}}

\newcommand{\bef}{\begin{figure}}
\newcommand{\enf}{\end{figure}}

\title{A New Model for Solving Narrow Escape Problem in Domain with Long Neck\thanks{This
        work was supported by Korean National Research Foundation, No. 2009-0085987 and BK21+ at Inha University.}}

\author{
Xiaofei Li\thanks{\footnotesize Department of Mathematics, Inha University, Incheon
402-751, Korea ( {\tt xiaofeili@inha.edu}).}}

\begin{document}

\maketitle

\begin{abstract}
The narrow escape problem arises in deriving the asymptotic expansion of the
solution of an inhomogeneous mixed Dirichlet-Neumann boundary value problem. In this paper, we mainly deal with narrow escape problem in a smooth domain connected to a long neck-Dendritic spine shape domain, which has a certain significance in biology. Since the special geometry of dendritic spine, we develop a new model for solving this narrow escape problem which is Neumann-Robin Boundary Model. This model transform spine singular domain to smooth spine head domain by inserting Robin boundary condition to the connection part between spine head and neck. We rigorously find the
high-order asymptotic expansion of Neumann-Robin Boundary Model and apply it to the solution of narrow escape problem in a dendritic spine shape domain. Our results show that the asymptotic expansion of the Neumann-Robin Boundary Model can be easily applied to the narrow escape problem for any smooth spine head domain with straight spine neck. By numerical simulations, we show that there is great agreement between the results of our Neumann-Robin Boundary Model and the original escape problem. In this paper, we also get some results for non-straight long spine neck case by considering curvature of spine neck.
\end{abstract}



\pagestyle{myheadings}
\thispagestyle{plain}
\markboth{XIAOFEI LI}{NEUMANN-ROBIN BOUNDARY MODEL}

\section{Introduction}
When a Brownian particle is confined in a bounded domain with a small absorbing windows on an otherwise reflecting boundary, it attempts to escape from this domain through this small absorbing windows. Narrow escape problem is to calculate the mean first passage time Brownian particle takes to get to the absorbing window. From the biological point of view,
the Brownian particles could be diffusing ions, globular proteins or cell-surface receptors. It is then of interest to
determine, for example, the mean time that an ion requires to find an open ion channel located in the cell membrane
or the mean time of a receptor to hit a certain target binding site.

In two dimension, the results of narrow escape problem in smooth bounded domain with one absorbing window was relatively complete \cite{acklm1,acklm4,acklm9,acklm15,acklm17,acklm19,acklm22,acklm23}. When there are several absorbing windows on the boundary, interaction of multiple absorbing windows are discussed in \cite{acklm10,acklm12}. In paper \cite{acklm21}, several kinds of singular domains have been discussed. In three dimension, the case that bounded domain is a ball with spherical boundary has been discussed in \cite{acklm5,acklm7}.

 While, our interest is different from these talked above, but another example of singular domain that a smooth domain connected by a long neck, such as dendritic spine(Fig.\ref{spine}). Dendritic spines serve as a storage site for synaptic strength and help transmit electrical signals to the neuron's cell body. As the important site of excitatory synaptic interaction, dendritic spines play an important role in neural plasticity, and their ability to regulate calcium attracts interests of many mathematicians and biologists \cite{acklm2,acklm3,acklm6,acklm8,acklm14}. Each spine has a bulbous head, and a thin neck that connects the head of the spine to the shaft of the dendrite. We consider simplified model of calcium diffusion in dendritic spines, which is discussed in \cite{acklm13}. That is, first we consider the calcium irons to be point charges, furthermore, we assume the motion of irons is free Brownian motion; second, the interaction between two electrostatic ions is neglected; third, we shall simply ignore impenetrable obstacles to the ionic motion posed by the presence of proteins. Thus, the iron motion inside the dendritic spine is geometrically unrestricted. In this paper, we regard the iron as calcium molecule. The calcium diffusion problem is narrow escape problem that is approximated by free Brownian motion in a domain which consists of a spherical head whose length is $L$ and a long cylindrical neck whose radius is $a$, where the radius of the neck is sufficiently small relative to that of the spine head(Fig.\ref{spinedomain}).

In this paper, we only talk about two dimensional case, where spine head is a bounded domain with smooth boundary and spine neck is rectangle, while, Neumann-Robin Boundary Model(\ref{NeumannRobin}) can be easily applied to three dimensions.

The narrow escape problem can be explained explicitly in the following way. Let $\Omega$ be a bounded simply connected domain in $\mathbb{R}^{2}$. Suppose
that $\partial \Omega$ is decomposed into the reflecting part $\partial \Omega_{r}$ and the absorbing part $\partial \Omega_{a}$. We assume that $\varepsilon=|\partial \Omega_{a}|/2$  is much smaller than the whole  boundary(Fig.\ref{spinedomain}). The narrow escape problem is to calculate the mean first passage time which is the solution $u_{\varepsilon}$ to (\ref{main_equation}),
\bef
\centering
  \includegraphics[width=2in]{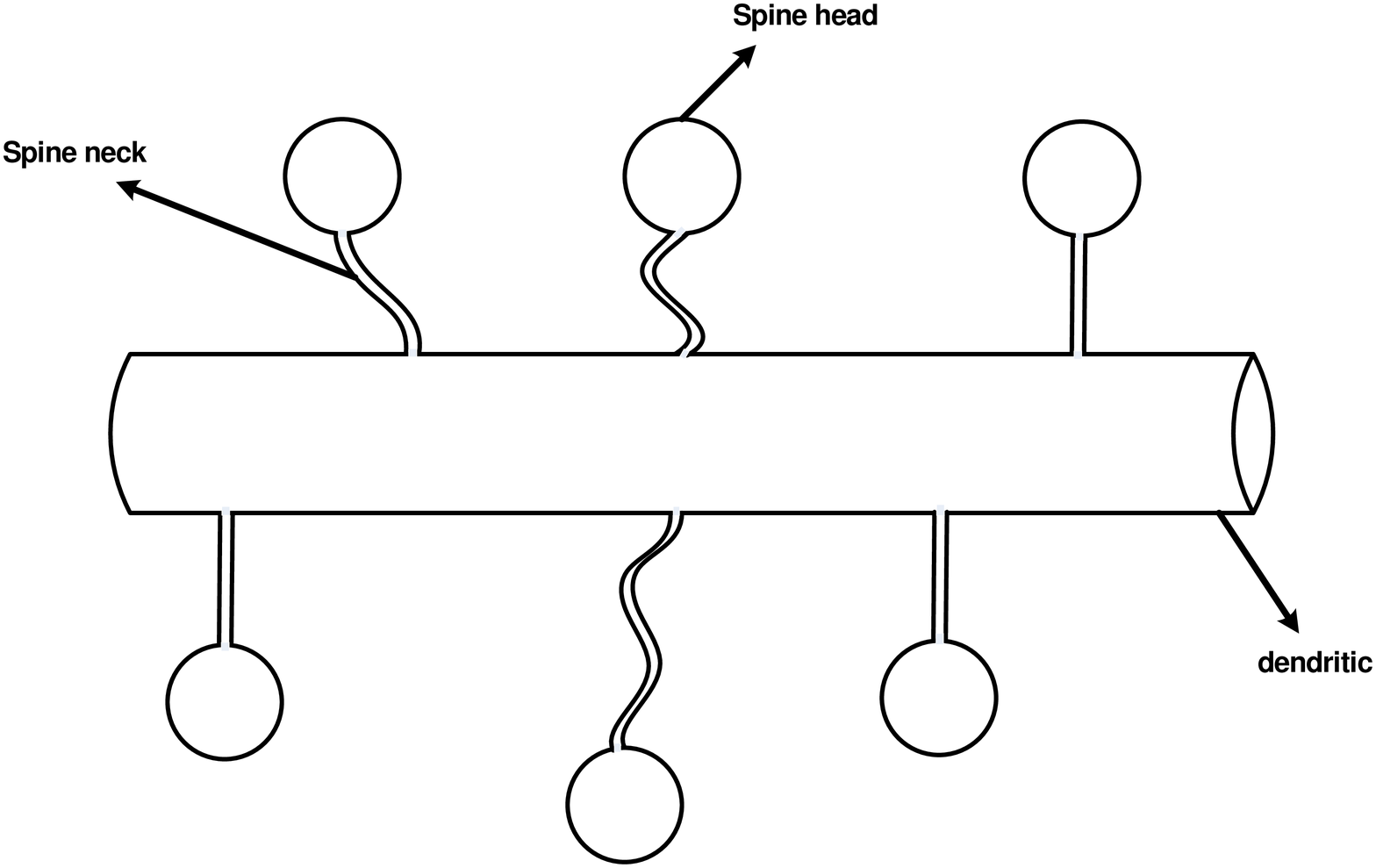}\\
  \caption{Abstract graph of dendritic spine}
  \label{spine}
 \enf

 \bef
\centering
  \includegraphics[width=1.5in]{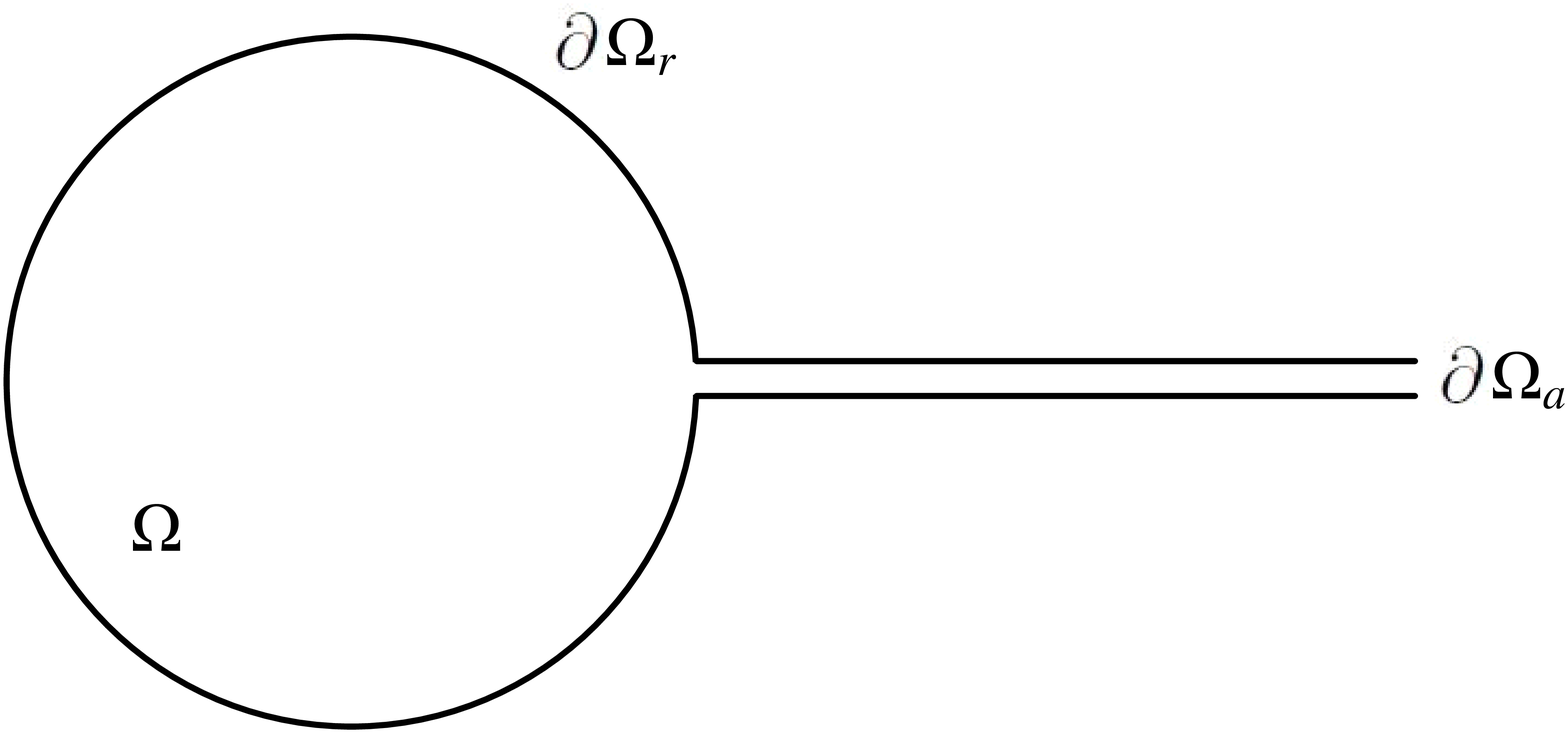}\\
  \caption{The modeling shape of dendritic spine with disk spine head and long spine neck, where $\Omega$ is the domain with long neck, $\partial \Omega_{r}$ is the reflection part, $\partial \Omega_{a}$ is the absorbing part.}
  \label{spinedomain}
 \enf

\beq \label{main_equation}
   \begin{cases}
         \triangle u_{\varepsilon}=-1,\quad &\mbox{in}~\Omega,\\
         \frac{\partial u_{\varepsilon}}{\partial \nu}=0,&\mbox{on}~\partial \Omega_{r},\\
         u_{\varepsilon}=0,&\mbox{on}~\partial \Omega_{a}.
         \end{cases}
         \eeq
The asymptotic analysis for narrow escape problem arises in deriving the asymptotic expansion of $u_{\varepsilon}$ as
$\varepsilon \rightarrow 0$, from which one can estimate the escape time of the Brownian particle.

In this work, instead of \eqnref{main_equation}, we develop another proper model to solve narrow escape problem in dendrite spine shape domain which we call it Neumann-Robin Boundary Model. The model is described by the following equations in domain $\Omega_h$(Fig.\ref{spinehead}),

\beq \label{NeumannRobin}
   \begin{cases}
         \triangle u_{\varepsilon}=-1,\quad &\mbox{in}~\Omega_h,\\
         \frac{\partial u_{\varepsilon}}{\partial \nu}=0,&\mbox{on}~\partial \Omega_{r},\\
         \frac{\partial u_{\varepsilon}}{\partial \nu}+\alpha u_{\varepsilon}=\beta,&\mbox{on}~\Gamma_\varepsilon:=\partial \Omega_h \setminus \overline{\partial \Omega_{r}}.
         \end{cases}
         \eeq
where $\Omega_h$ is the spine head of $\Omega$ mentioned in Fig.\ref{spine}, the size of $\Gamma_\varepsilon$ is still small with  $|\Gamma_\varepsilon|=2\varepsilon$, but it is not just an absorbing boundary any more. Note that, the domain $\Omega$ we considered here is only the head domain, without the neck.

 \bef
\centering
  \setlength{\unitlength}{1bp}%
  \begin{picture}(83.10, 81.02)(0,0)
  \put(0,0){\includegraphics{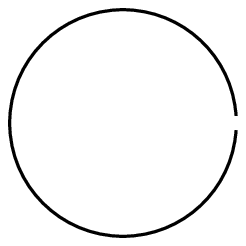}}
  \put(20.48,32.95){\fontsize{8.54}{10.24}\selectfont $\Omega_h$}
  \put(61.06,68.68){\fontsize{8.54}{10.24}\selectfont $\partial\Omega_r$}
  \put(73.17,32.34){\fontsize{8.54}{10.24}\selectfont $\Gamma_{\varepsilon}$}
  \end{picture}
    \caption{The domain considered in Neumann-Robin model, which omit the long neck by adding Robin boundary condition to small arc $\Gamma_\varepsilon$, $\partial \Omega_{r}$ still represents reflecting boundary.}
    \label{spinehead}
    \enf

In this paper, we analyze the asymptotic behavior of the solution of Neumann-Robin Boundary Model in the domain $\Omega_h$ (see Fig.\ref{spinehead}) in two dimensions. Actually, the asymptotic behavior can be applied to any smooth bounded domain in two dimensions, not only disk case. We rigorously derived the expansion formula for (1.2) up to order $O(\varepsilon)$
\beq
u_{\varepsilon}(x)=\frac{|\Omega|}{2\alpha\varepsilon}+\frac{|\Omega|}{\pi}\left(\frac{3}{2}+\ln\frac{1}{2\varepsilon}\right)+ \frac{\beta}{\alpha}+\Phi_{\Omega}(x,x^*)+O(\varepsilon).
\end{equation}
where $\Phi_{\Omega}(x,x^*)$ can be referred to (\ref{Phi_x*}).

By assigning specific $\alpha=1/L$, $\beta=L/2$ to our Neumann-Robin Boundary Model(see the reason for choice of $\alpha,\beta$ in section 4), the solution formula (1.3) can approximate the mean first passing time of narrow escape problem (1.1) in spine domain (Fig.\ref{spinedomain}) up to order $O(\varepsilon)$. The numerical results show a great agreement between them.

This paper is organized as follows. In Section 2, we review the Neumann function for Laplacian and introduce an integral operator for further calculations. In section 3 the asymptotic formula for the solution
to Neumann-Robin Boundary Model has been rigorously derived by using layer potential techniques. In section 4 and 5, we discuss how the Neumann-Robin Boundary Model corresponds to the original escape problem theoretically and numerically. The paper ends with a short conclusion.

\section{Preliminaries}

Let $N(x,z)$ be the Neumann function for $-\triangle$ in $\Omega$ corresponding to a Dirac mass at $z\in\Omega$. We assume $\partial \Omega$ is $C^2$ smooth. $N(x,z)$ is the solution to
\beq
  \begin{cases}
        \triangle_{x} N(x,z)=-\delta_{z}, &x\in \Omega,\\
      \displaystyle\frac{\partial N}{\partial \nu_{x}}=-\frac{1}{|\partial \Omega|},  & x\in \partial\Omega,\\
  \end{cases}
\eeq
For uniqueness, we assume $\int_{\partial \Omega}N(x,z)d\sigma(x)=0$.

If $z\in \Omega$, $N(x,z)$ can be written in the form
\beq N(x,z)=-\frac{1}{2\pi}\ln|x-z|+R_{\Omega}(x,z),\quad x\in\Omega,\eeq
where $R_{\Omega}(x,z)$ is the regular part which belongs to $H^{3/2}(\Omega)$, and solves
\beq
\begin{cases}
-\triangle_{x} R_{\Omega}(x,z)=0,&x\in\Omega,\\
\displaystyle\frac{\partial R_{\Omega}}{\partial \nu_{x}}\Big|_{x\in \partial \Omega}=-\frac{1}{|\partial \Omega|}+\frac{1}{2\pi}\frac{\langle x-z,\nu_{x}\rangle}{|x-z|^{2}},&x\in \partial \Omega.
\end{cases}
\end{equation}\\
If $z\in \partial \Omega$, Neumann function on the boundary $N_{\partial \Omega}$, can be written as
\beq \label{N_p} N_{\partial \Omega}(x,z)=-\frac{1}{\pi}\ln|x-z|+R_{\partial \Omega}(x,z),\quad x \in \Omega,z \in \partial \Omega, \eeq
where the singularity of $N_\partial \Omega(x,z)$ is $-\frac{1}{\pi}\ln|x-z|$ (See \cite{acklm1}), $R_{\partial \Omega}(x,z)$  solves the problem
\beq
  \begin{cases}
        \triangle_{x}R_{\partial \Omega}(x,z)=0,& x\in \Omega,\\
       \displaystyle\frac{\partial R_{\Omega}}{\partial \nu_{x}}\Big|_{x\in \partial \Omega}=-\frac{1}{|\partial \Omega|}+\frac{1}{\pi}\frac{\langle x-z,\nu_{x}\rangle}{|x-z|^{2}},&x\in \partial \Omega, ~z\in \partial \Omega.\\
        \end{cases}
        \end{equation}
Note that the Neumann data above is bounded on $\partial \Omega$ uniformly in $z\in \partial \Omega$ since $\partial \Omega$ is $C^{2}$-smooth, and hence $R_{\partial \Omega}(x,z)$ belongs to $H^{3/2}(\Omega)$ uniformly in $z\in \partial \Omega$.(See \cite{acklm1})

For later use, we introduce the integral operator $L:L^2[-1,1]\rightarrow L^2[-1,1]$, defined by
$$L[\phi](x)=\int\limits_{-1}^{1}\ln|x-y|\phi(y)dy.$$
We can see operator $L$ is bounded (see Lemma 2.1 in \cite{acklm1}).

\section{Neumann-Robin boundary value problem}
In this section we rigorously give the asymptotic analysis of our Neumann-Robin Boundary Model in a more general spine head domain $\Omega_h$(Fig. \ref{gene}), where $\Omega_h\in C^2(R^2)$, and derive full expansion solution formula for (\ref{NeumannRobin}) in domain $\Omega_h$ up to order $O(\varepsilon)$.
  \bef
\centering
    \setlength{\unitlength}{1bp}%
  \begin{picture}(113.41, 75.37)(0,0)
  \put(0,0){\includegraphics{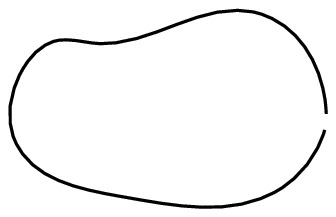}}
  \put(29.43,26.88){\fontsize{8.54}{10.24}\selectfont $\Omega_h$}
  \put(85.96,63.04){\fontsize{8.54}{10.24}\selectfont $\partial\Omega_{r}$}
  \put(99.72,25.35){\fontsize{8.54}{10.24}\selectfont $\Gamma_{\varepsilon}$}
  \end{picture}
    \caption{Any smooth domain with small opening $\Gamma_\varepsilon$. This domain is where Neumann-Robin Boundary Model considered.}
    \label{gene}
    \enf

 We consider the Laplace equation in $\Omega_h$ with the mixed Neumann-Robin boundary condition. The Robin boundary condition is imposed on $\Gamma_{\varepsilon}$($\Gamma_{\varepsilon}$ is a very small part) and the Neumann boundary condition on the part $\partial \Omega_{r}:= \partial \Omega_h \setminus \overline{\Gamma_\varepsilon}$  :
 \beq  \label{NR_eq}
   \begin{cases}
         \triangle u_{\varepsilon}=-1, &\mbox{in}~\Omega_h,\\
         \frac{\partial u_{\varepsilon}}{\partial \nu}=0,&\mbox{on}~\partial \Omega_{r},\\
         \frac{\partial u_{\varepsilon}}{\partial \nu}+\alpha u_{\varepsilon}=\beta,&\mbox{on}~\Gamma_{\varepsilon}.
         \end{cases}
         \end{equation}
Here, $\alpha>0$ and $\beta$ are given constants. We assume that $\varepsilon$ is sufficiently small so that $\alpha \varepsilon \ll 1$ and $\alpha< \alpha_0$ for some constant $\alpha_0>0$.

The goal in this section is to derive the asymptotic expansion of $u_{\varepsilon}$ as $\varepsilon\rightarrow 0$, from which one can estimate the exit time of the calcium iron in the spine head.

By integrating the first equation in (3.1) over $\Omega_h$ using the divergence theorem we get
\beq
\int_{\Gamma_{\varepsilon}}\frac{\partial u_{\varepsilon}}{\partial \nu}d\sigma=-|\Omega_h|.
\end{equation}

Let us define $g(x)$ by
$$g(x)=\int_{\Omega_h}N(x,z)dz,\quad x\in\Omega_h,$$\\
which satisfies
\beq \label{g_eq}
\begin{cases}
\triangle g=-1,& \mbox{in}~\Omega_h,\\
\displaystyle\frac{\partial g}{\partial\nu}=-\frac{|\Omega_h|}{|\partial \Omega_h|}, & \mbox{on}~ \partial \Omega_h,\\
\displaystyle \int_{\partial \Omega_h} g d\sigma=0.&
\end{cases}
\end{equation}
Therefore, applying the Green's formula  to  $u_\varepsilon$ and the Neumann function $N$ and using \eqnref{NR_eq} and \eqnref{g_eq}, we get
\beq \label{u_form1}
u_{\varepsilon}(x)=g(x) + \int_{\Gamma_{\varepsilon}}N_{\partial\Omega_{h}}(x,z)\frac{\partial u_{\varepsilon}(z)}{\partial \nu_{z}}d\sigma(z)
+C_{\varepsilon},
\eeq
where $$C_{\varepsilon}=\frac{1}{|\partial \Omega_h|}\int_{\partial \Omega_h}u_{\varepsilon}(z)d\sigma(z).$$
By \eqnref{N_p}, the equation \eqnref{u_form1} becomes
\beq
u_{\varepsilon}(x)=g(x)-\frac{1}{\pi}\int_{\Gamma_{\varepsilon}}\ln|x-z|\frac{\partial u_{\varepsilon}(z)}{\partial \nu}d\sigma(z)
+ \int_{\Gamma_{\varepsilon}}R_{\partial \Omega_h}(x,z)\frac{\partial u_{\varepsilon}(z)}{\partial \nu}d\sigma(z)+C_{\varepsilon}.\\
 \end{equation}
On $\Gamma_\varepsilon$, by Robin boundary condition, (3.5) can be written as
     \beq \label{Robin_bc}
        \frac{1}{\pi}\int_{\Gamma_{\varepsilon}}\ln|x-z|\frac{\partial u_{\varepsilon}(z)}{\partial \nu_{z}}d\sigma(z)
        - \int_{\Gamma_{\varepsilon}}R_{\partial \Omega_h}(x,z)\frac{\partial u_{\varepsilon}(z)}{\partial \nu_{z}}d\sigma(z)
        =\frac{1}{\alpha}\frac{\partial u_{\varepsilon}(x)}{\partial \nu_{x}}-\frac{\beta}{\alpha}+g(x)+C_{\varepsilon}.
     \eeq

Let $x(t):[-\varepsilon,\varepsilon]\rightarrow \mathbb{R}^{2}$  be the arc-length parametrization of $\Gamma_{\varepsilon}$, {\it i.e.},
         $|x'(t)|=1$ for all $t \in [-\varepsilon,\varepsilon]$ and
        $$\Gamma_{\varepsilon}=\{\,x(t)~|~t\in [-\varepsilon,\varepsilon]\,\}.$$
For simplicity, we let
\beq f(t)=g(x(t)), \quad\phi_{\varepsilon}(t)=\frac{\partial u_{\varepsilon}}{\partial \nu}(x(t)), \quad r(t,s)=R_{\partial \Omega_h}(x(t),x(s)).\eeq
Then it follows from \eqnref{Robin_bc} that
 \beq
\frac{1}{\pi}\int_{-\varepsilon}^{\varepsilon}\ln|\,x(t)-x(s)|\phi_{\varepsilon}(s)ds
- \int_{-\varepsilon}^{\varepsilon}r(t,s)\phi_{\varepsilon}(s)ds
=\frac{1}{\alpha}\phi_{\varepsilon}(t)+f(t)+C_{\varepsilon}-\frac{\beta}{\alpha}.
        \end{equation}

By the change of variable, we obtain
\beq \label{main_int_eq}
\frac{1}{\alpha \varepsilon}\widetilde{\phi}_{\varepsilon}(t)-\frac{1}{\pi}\int_{-1}^{1}\ln|x(\varepsilon t)-x(\varepsilon s)|\widetilde{\phi}_{\varepsilon}(s)ds
+ \int_{-1}^{1}r(\varepsilon t,\varepsilon s)\widetilde{\phi}_{\varepsilon}(s)ds
=-f(\varepsilon t)-C_{\varepsilon}+\frac{\beta}{\alpha},
        \end{equation}
where $\widetilde{\phi}_{\varepsilon}(t)=\varepsilon \phi_{\varepsilon}(\varepsilon t)$.

We define two bounded integral operators $L,~ L_1: L^2[-1,1]\rightarrow L^2[-1,1]$ by
\begin{align*} L[\phi]&=\int_{-1}^1 \ln|t-s|\phi(s)ds,\\
L_{1}[\phi]&=\frac{1}{\varepsilon}\int_{-1}^{1}\left(\ln\frac{|x(\varepsilon t)-x(\varepsilon s)|}{\varepsilon |t-s|}+\pi r(0,0)-\pi r(\varepsilon t,\varepsilon s)\right)\phi(s)ds.
\end{align*}
Since
$|x(\varepsilon t)-x(\varepsilon s)|=\varepsilon|t-s|(1+O(\varepsilon)),$ one can see that $L_{1}$ is bounded independently of $\varepsilon$.

Using the compatibility condition
\beq \label{compat_cond} \int_{-1}^{1}\widetilde{\phi}_{\varepsilon}(t)dt=-|\Omega_h|,\eeq
we may write \eqnref{main_int_eq} as
\beq
\frac{1}{\alpha \varepsilon}\widetilde{\phi}_{\varepsilon}(t)-\frac{1}{\pi}(L+\varepsilon L_{1})[\widetilde{\phi}_{\varepsilon}](t)=-\frac{|\Omega_h|\ln\varepsilon}{\pi}+r(0,0)|\Omega_h|-f(\varepsilon t)-C_{\varepsilon}+\frac{\beta}{\alpha}.
\eeq

Assume $\alpha< \alpha_0$ and $\alpha \varepsilon \ll 1$. Then we have
\beq
(I-\frac{\alpha\varepsilon}{\pi}(L+\varepsilon L_{1}))^{-1}=I+\frac{\alpha\varepsilon}{\pi}(L+\varepsilon L_{1})+O(\alpha^2\varepsilon^{2}).
\eeq
Therefore we have
\begin{align}
\widetilde{\phi}_{\varepsilon}(t)
&={-\alpha\varepsilon}\left[I+\frac{\alpha\varepsilon}{\pi}(L+\varepsilon L_{1})+O(\alpha^2\varepsilon^{2})\right]\left(C_{\varepsilon}+\frac{|\Omega_h|\ln\varepsilon}{\pi}- r(0,0)|\Omega_h|+ f(\varepsilon t)-\frac{\beta}{\alpha}\right)\\
&={-\alpha\varepsilon}\left[I+\frac{\alpha\varepsilon}{\pi}(L+\varepsilon L_{1})+O(\alpha^2\varepsilon^{2})\right]\left(\tilde{C}_{\varepsilon}+O(\varepsilon)\right),
\end{align}
where $\tilde{C}_\varepsilon=C_{\varepsilon}+\frac{|\Omega_h|\ln\varepsilon}{\pi}- r(0,0)|\Omega_h|+ f(0)-\frac{\beta}{\alpha}$.

 By \eqnref{compat_cond}, we see $  \tilde{C}_\varepsilon = O((\alpha \varepsilon)^{-1})$. Then collecting terms we have
\beq \label{tildephi}
\widetilde{\phi}_{\varepsilon}(t)
=-{\alpha\varepsilon}\tilde{C}_\varepsilon -\frac{(\alpha\varepsilon)^2}{\pi}\tilde{C}_\varepsilon L[1](t)+O(\alpha\varepsilon^{2}).
\eeq
Plugging it into \eqnref{compat_cond} we obtain
\beq
2{\alpha\varepsilon}\tilde{C}_\varepsilon +\frac{(\alpha\varepsilon)^2}{\pi}\tilde{C}_\varepsilon\int_{-1}^1 L[1](t)dt = {|\Omega_h|}+O(\alpha\varepsilon^{2}).
\eeq
Then we get
\begin{align*}
\tilde{C}_\varepsilon &= \left(1 +\frac{\alpha\varepsilon}{2\pi}\int_{-1}^1 L[1](t)dt \right)^{-1}\left(\frac{|\Omega_h|}{2\alpha\varepsilon }+O(\varepsilon)\right)\\
&=\frac{|\Omega_h|}{2\alpha\varepsilon }-\frac{|\Omega_h|}{4\pi}\int_{-1}^1 L[1](t)dt+O(\varepsilon).
\end{align*}

The direct calculation shows us
$$\int_{-1}^{1}L[1](t)dt=\int_{-1}^{1}\int_{-1}^{1}\ln|t-y|dtdy=4\ln2-6.$$
Therefore we arrive at
\beq \label{tildeC}
\tilde{C}_\varepsilon=\frac{|\Omega_h|}{2\alpha\varepsilon}+\frac{|\Omega_h|}{\pi}(\frac{3}{2}-\ln2)+O(\varepsilon),
\end{equation}
and hence
\beq \label{C_exp}
C_{\varepsilon}=\frac{|\Omega_h|}{2\alpha\varepsilon}+\frac{|\Omega_h|}{\pi}(\frac{3}{2}+\ln\frac{1}{2\varepsilon})+ \frac{\beta}{\alpha}+r(0,0)|\Omega_h|-f(0)+O(\varepsilon).
\end{equation}

Substituting \eqnref{tildeC} into \eqnref{tildephi}, we have
\beq
\widetilde{\phi}_{\varepsilon}(t)=-\frac{|\Omega_h|}{2}-\frac{|\Omega_h|}{\pi}\alpha\varepsilon \left(\frac{3}{2}-\ln2+\frac{1}{2}L[1](t)\right)+O(\alpha\varepsilon^{2}),
\eeq
 and
\beq
\phi_{\varepsilon}=\frac{1}{\varepsilon}\widetilde{\phi}_{\varepsilon}(\frac{t}{\varepsilon})
=-\frac{|\Omega|}{2\varepsilon}-\frac{|\Omega|}{\pi}\alpha\left(\frac{3}{2}-\ln2+\frac{1}{2}L[1]\left(\frac{t}{\varepsilon}\right)\right)+O(\alpha\varepsilon^{3/2}),
\end{equation}
where $O(\alpha\varepsilon^{2})$ and $O(\alpha\varepsilon^{3/2})$ are measured in $\parallel \cdot \parallel_{L^2[-1,1]}$ and $\parallel \cdot \parallel_{L^2[-\varepsilon,\varepsilon]}$, respectively.

If $x$ is away from $\Gamma_{\varepsilon}$, {\it i.e.}, $\mbox{dist}(x,\Gamma_{\varepsilon})\geq c$ for some constant $c>0$, then
\begin{align}
&\int_{\Gamma_{\varepsilon}}N_{\partial \Omega_h}(x,z)\frac{\partial u_{\varepsilon}(z)}{\partial \nu_{z}}d\sigma(z)\nonumber\\
&=\int_{-\varepsilon}^\varepsilon N_{\partial \Omega_h}(x,z(t))\left(-\frac{|\Omega_h|}{2\varepsilon}-\frac{|\Omega_h|}{\pi}\alpha \left(\frac{3}{2}-\ln2+\frac{1}{2}L[1]\left(\frac{t}{\varepsilon}\right)\right)+O(\alpha\varepsilon^{3/2})\right)dt\nonumber\\
&=-|\Omega_h| N_{\partial\Omega_h}(x,x^*)+O(\epsilon).\label{phi_exp}
\end{align}

Finally, combining \eqnref{u_form1}, \eqnref{C_exp} and \eqnref{phi_exp} yields
\beq
  \begin{split}
u_{\varepsilon}(x) &=g(x)+\int_{\Gamma_{\varepsilon}}N_{\partial \Omega_h}(x,z)\frac{\partial u_{\varepsilon}(z)}{\partial \nu_{z}}d\sigma(z)+C_{\varepsilon}\\
 &=g(x)-|\Omega|N_{\partial \Omega_h}(x,x^*)+\frac{|\Omega_h|}{2\alpha\varepsilon}+\frac{|\Omega_h|}{\pi}(\frac{3}{2}+\ln\frac{1}{2\varepsilon})+ \frac{\beta}{\alpha}+r(0,0)|\Omega_h|-f(0)+O(\varepsilon)
\end{split}
\end{equation}
for $x\in \Omega_h$ provided that dist(x,$\Gamma_{\varepsilon}$)$\geq c$ for some constant $c>0$. Thus we have the following theorem.\\
\begin{thm} Suppose that $\Gamma_{\varepsilon}$ is an arc of center $x^{*}$ and length $2\epsilon$. Then the following asymptotic expansion of $u_{\epsilon}$ for (3.1) holds
\beq
u_{\varepsilon}(x)=\frac{|\Omega_h|}{2\alpha\varepsilon}+\frac{|\Omega_h|}{\pi}\left(\frac{3}{2}+\ln\frac{1}{2\varepsilon}\right)+ \frac{\beta}{\alpha}+\Phi_{\Omega_h}(x,x^*)+O(\varepsilon),
\end{equation}
where
\beq \label{Phi_x*}
\Phi_{\Omega_h}(x,x^*)=\int_{\Omega_h}N(x,z)dz-|\Omega_h|N_{\partial \Omega_h}(x,x^*)-\int_{\Omega_h}N(x^*,z)dz+|\Omega_h|R_{\partial \Omega_h}(x^*,x^*).
\end{equation}
The remainder $O(\varepsilon)$ is uniform in $x\in \Omega_h$ satisfying dist$(x,\Gamma_{\varepsilon})\geq c$ for some constant $c>0$. Moreover, if $x(t),-\epsilon<t<\epsilon$, is the arclength parameterization of $\Gamma_{\varepsilon}$, then,
\beq
\frac{\partial u_{\epsilon}}{\partial \nu}(x(t))=-\frac{|\Omega_h|}{2\varepsilon}-\frac{|\Omega_h|}{\pi}\alpha\left(\frac{3}{2}-\ln2+\frac{1}{2}L[1](\frac{t}{\varepsilon})\right)+O(\alpha\varepsilon^{3/2}),\\
\end{equation}
where $O(\alpha\varepsilon^{3/2})$ is with respect to $\parallel$ $\parallel_{L^2[-\varepsilon,\varepsilon]}$.
\end{thm}

We note that the  function $\Phi_{\Omega_h}(x,x^*)$ solves the following problem
\beq
  \begin{cases}
        $$\triangle_{x}\Phi_{\Omega_h}(x,x^*) =0, &x\in \Omega_h,$$\\
        $$\frac{\partial \Phi_{\Omega_h}(x,x^*)}{\partial \nu_{x}}=-|\Omega_h|\delta_{x^{*}}, &x\in \partial \Omega_h.$$\\
        \end{cases}
        \end{equation}
If $\Omega_h$ is a unit disk centered at 0, one can easily see from (2.6) and (2.7) that
$$\Phi_{\Omega_h}(x,x^{*})=\ln|x-x^{*}|+\frac{1}{4}(1-|x|^{2}).$$

\section{Calcium diffusion in dendritic spines}
In this section we use our new Neumann-Robin Boundary Model to solve the narrow escape problem in dendritic spine shape domain(Fig. \ref{spinedomain}) which is calcium diffusion problem. That is, we calculate how long a single calcium molecule stays in the spine before it escapes from it.

The usual way to calculate the solution in smooth domain requires boundary layer expansions for small window size and the asymptotic of the Neumann function that worked for nonsingular problems failed for calcium diffusion model since there are two singular points on the conneting part of spine head and spine neck. A quite different approach to the asymptotic problem is required which are much different from those reported in the cited reviews.

In order to approximate the escape time for a particle in spine head, we approach a new method, different from those dealt with in \cite{acklm11,acklm16}, that we use the Neumann-Robin Boundary Model in the spine head domain(Fig.\ref{spinehead}), but with the specific
$$\alpha=\frac{1}{L},\beta=\frac{L}{2},$$
on the boundary $\Gamma_{\epsilon}$ which is the opening part of the big head. Here $L$ is the length of the neck. Note that we have changed the domain by dropping the long neck and assigning Robin boundary condition to the connecting arc $\Gamma_{\varepsilon}$ between the spine head and the long neck. Instead of dealing with the singular part on spine domain, we put a proper Robin boundary to the connecting part between spine head and neck.
\bef
\centering
    \includegraphics[width=2.6in]{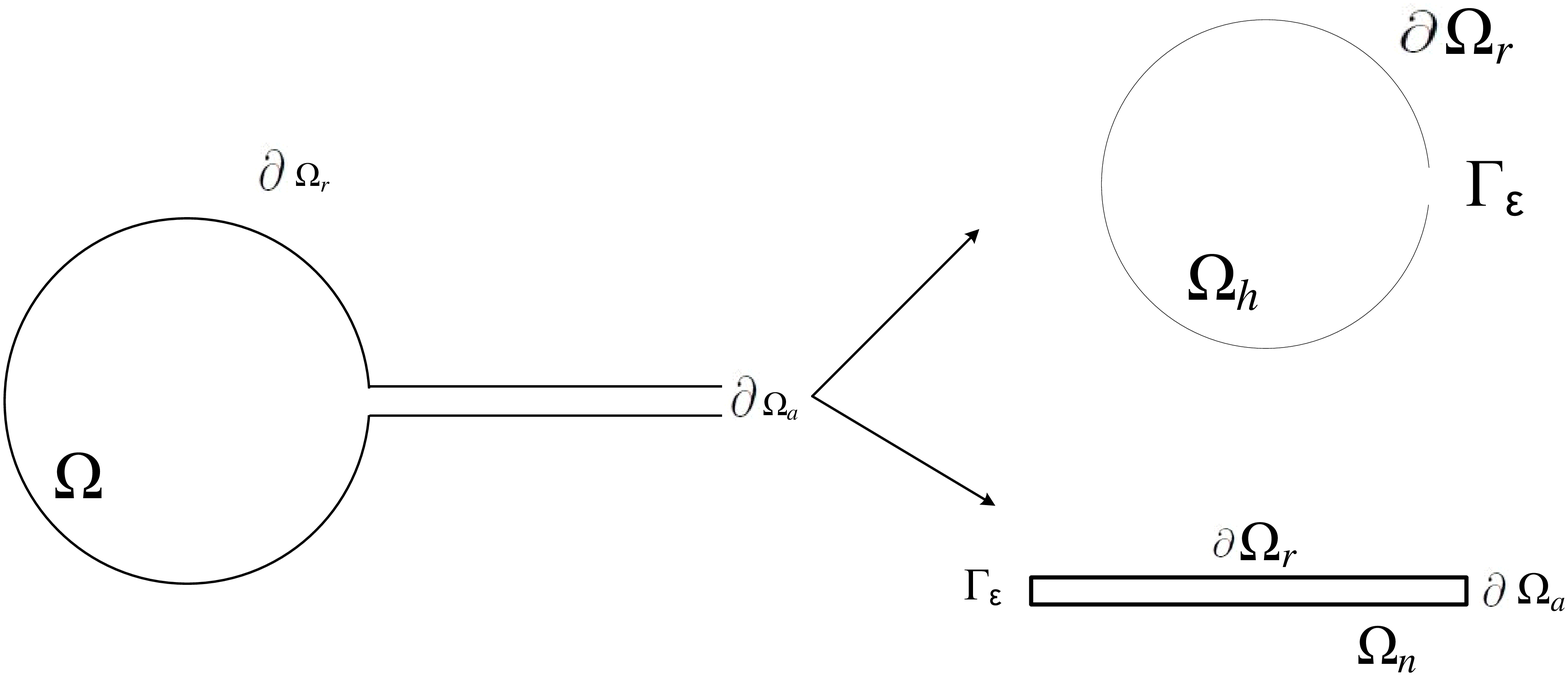}\\
    \caption{Decompose spine domain into two parts $\Omega=\Omega_{h}+\Omega_{n}$. One is spine head domain $\Omega_{h}$, with reflecting boundary $\partial \Omega_r$ and Robin boundary $\Gamma_\varepsilon$ (where the Neumann-Robin Model can be solved using layer potential techniques in this smooth domain), the other part is long neck $\Omega_{n}$.}
    \label{decompose}
    \enf

The heuristic reason for this specific choice of $\alpha$ and $\beta$ comes from the following:

In the spine domain, we decompose the domain into two parts, one is spine head which has smooth boundary, the other is the long neck domain(Fig.\ref{decompose}). Since the spine neck radius is small enough,  we assume the escape time on the small part $\Gamma_{\varepsilon}$(which connects head and long neck) is constant. Thus, in the spine neck domain $\Omega_n$(Fig.\ref{decompose}), escape time $u_\varepsilon$ satisfies the following equation
        \beq
        \begin{cases}
              \triangle u_\varepsilon=-1, &\mbox{in}~\Omega_{n},\\
              \frac{\partial u_\varepsilon}{\partial \nu} = 0, &\mbox{on}~\partial \Omega_{r},\\
              u_\varepsilon=0,  &\mbox{on}~\partial \Omega_{a},\\
              u_\varepsilon=C,  &\mbox{on}~\Gamma_{\varepsilon},
        \end{cases}
      \end{equation}
where $C$ is constant whose value means the escape time for the point initiated on $\Gamma_{\varepsilon}$. Taking the center point of $\Gamma_{\varepsilon}$ to be original point $(0,0)$, by separation of variables, we can solve this partial differential equation in $\Omega_{n}$. The solution of (4.1) is
\beq
u_\varepsilon(x,y)=-\frac{1}{2}(L-x)^2+(\frac{C}{L}+\frac{L}{2})(L-x),
\eeq
where $x\in[0,L]$, $y\in(-\varepsilon,\varepsilon)$.
Since $\Omega_h$ and $\Omega_n$ are connected by $\Gamma_\varepsilon$, they share the same boundary value. The above solution $u_\varepsilon(x,y)$ satisfies the Robin boundary condition
\beq
\frac{\partial u_{\varepsilon}}{\partial \nu}+\alpha u_{\varepsilon}=\beta,
\eeq
with $$\alpha=\frac{1}{L},\beta=\frac{L}{2}.$$

Then the approximated Neumann-Robin Boundary Model for narrow escape problem is:
\beq
   \begin{cases}
         \triangle u_{\varepsilon}=-1, &\mbox{in}~\Omega_{h},\\
         \frac{\partial u_{\varepsilon}}{\partial \nu}=0, &\mbox{on}~\partial \Omega_{r},\\
         \frac{\partial u_{\varepsilon}}{\partial \nu}+\frac{u_{\varepsilon}}{L}=\frac{L}{2}, &\mbox{on}~\Gamma_{\varepsilon},
         \end{cases}
         \end{equation}
where $u_{\varepsilon}(x)$ is the escape time of the calcium iron which initiated at $x$ position in the spine head. $\Omega_{h}$ is the domain of the dendritic spine head and $\partial \Omega_{r}$ means the boundary where calcium molecule is reflected.
According to Theorem 3.1 in the last section, the solution to (4.2) is:
\beq
u_{\epsilon}(x)=\frac{|\Omega_{h}|L}{2\varepsilon}+\frac{|\Omega_{h}|}{\pi}(\frac{3}{2}+\ln\frac{1}{2\varepsilon})+\frac{L^2}{2}+\Phi_{\Omega_{h}}(x,x^*)+O(\varepsilon),
\end{equation}
where $\Phi_{\Omega_{h}}(x,x^*)$ is the same as (3.24) for the domain $\Omega_h$.

Suppose that $\Omega$ is the spine with straight spine neck domain. The length of the neck is $L$, and $\partial\Omega_a$ is the exiting arc(See Fig.\ref{spinedomain}). $\Omega_{h}$ is the spine head, and $\Gamma_{\varepsilon}$ is the arc of center $x^*$ which connects spine head and spine neck. The first mean passage time $u_\varepsilon(x)$ of a Brownian particle confined in $\Omega_h$ exiting through $\partial\Omega_a$ can be approximated by the following formula
\beq
u_\varepsilon(x)\approx\frac{|\Omega_{h}|L}{2\varepsilon}+\frac{|\Omega_{h}|}{\pi}(\frac{3}{2}+\ln\frac{1}{2\varepsilon})+ \frac{L^{2}}{2}+\Phi_{\Omega_{h}}(x,x^*),
\end{equation}
where $\Phi_{\Omega_{h}}(x,x^*)$ is the same as (3.24) for the domain $\Omega_h$. The error between formula (4.6) and the exact solution to (1.1) in spine domain is of order $O(\varepsilon)$, which can be seen by the numerical experiment data in the next section.

The method using Neumann-Robin Boundary Model to solve narrow escape problem in domain with long neck is quite different from what has been discussed in \cite{acklm11}. Their idea is to calculate the exit time by separating the exiting process of the particle into two processes. One is the time from the head to the interface $\Gamma_{\varepsilon}$ between head and neck, the other is the time from the interface to the absorbing arc. The mean first passage time can be obtained by adding the time of these two processes together. Their approximated formulation for planar spine connected to the neck at a right angle is
\beq
u_{\varepsilon}(x)=\frac{|\Omega_h|}{\pi}\ln\frac{|\partial\Omega_h|}{2\varepsilon}+O(1)+\frac{L^2}{2}+\frac{|\Omega_h|L}{2\varepsilon},
\eeq
where $O(1)$ is the error term. From (4.6) and (4.7), we can see that the results from these two methods have the same first leading order term $\frac{L^2}{2}+\frac{|\Omega_1|L}{2\varepsilon}+\frac{|\Omega_h|}{\pi}\ln\frac{1}{2\varepsilon}$. Thus, using our method we can obtain the exact formula for $O(1)$ in (4.7).
\section{Numerical experiment}
In order to check whether the asymptotic formula (4.6) can solve the narrow escape problem, we compare the numerical results of (4.6) denoted by $u_\varepsilon$ in spine head domain, with the numerical solutions obtained by solving the two dimensional narrow escape problem (1.1) by using Matlab, which is denoted in this section by $u$. Without loss of generality, we use the approximated spine geometry with unit disk spine head, and a rectangle neck in two dimension.
\bef
\centering
    \includegraphics[width=6in]{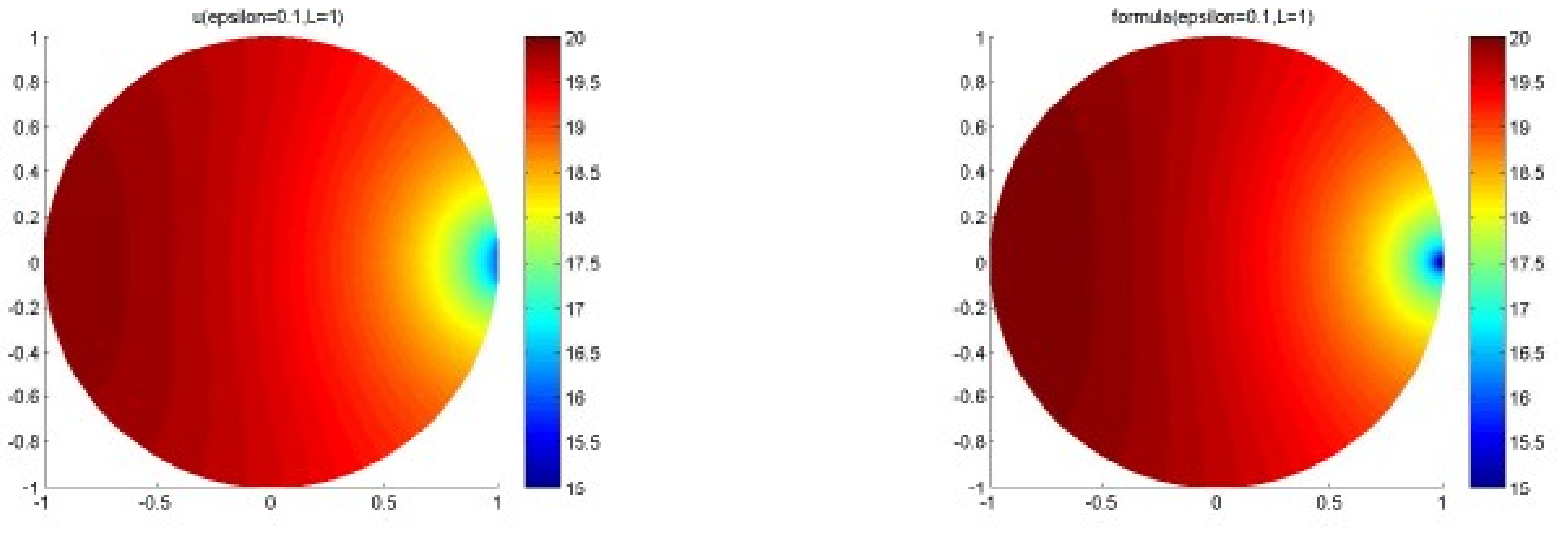}\\
    \caption{The spine geometry is approximated by a unit disk head and a rectangle neck with width $2\varepsilon=0.2$ and length $L=1$. Coordinates represent the position of the iron where it is initiated. The color represents the exit time of particle initiated at this point. Left figure: the numerical result $u$ of (1.1). Right figure:  exit time $u_\varepsilon$ computed by asymptotic formula (4.6).}
    \label{com01}
    \enf

The first goal in this section is to compare the expansion formula (4.6) and the numerical solution to (1.1). The domain of the escape problem (1.1) is given in Fig.\ref{spinedomain}. As one example, we choose the spine head to be unit disk, with the neck length $L=1$, and with the exit arc length $|\partial\Omega_a|=2\varepsilon$, $\varepsilon=0.1$. Then, the Neumann-Robin model solve the narrow escape problem in the spine head domain(Fig.\ref{spinehead}). Instead of considering the neck, we put the Robin boundary condition $\frac{\partial u_{\varepsilon}}{\partial \nu}+\frac{u_{\varepsilon}}{L}=\frac{L}{2}$ on $\Gamma_\varepsilon$, where the arc length of $\Gamma_\varepsilon$ is $2\varepsilon$, $\varepsilon=0.1$. Note that these two problems have the same spine head domain.

The numerical results are given in Fig. \ref{com01}. The figure on the left side gives the numerical solution to (1.1) in the spine head domain. The value at each point shows the exit time of the particle initiated at this point. We can easily see that if the particle is initiated near the small arc, then it takes less time to escape. On the other hand, the figure on the right side is the numerical result of the expansion formula (4.6). The value at every point represents the exit time of the particle initiated at this point. From these two figures, we can easily see the results in both situations coincide perfectly. The numerical data show that the error between these two situations are of order $O(\varepsilon)$. This will be seen in Table.\ref{tableu}.
\bef
\centering
    \includegraphics[width=2.2in]{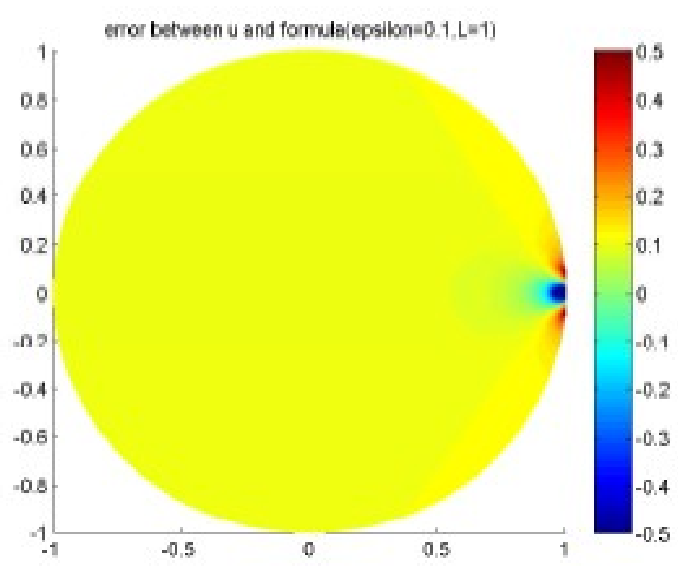} \\
    \caption{Error between $u$ and $u_\varepsilon$ in the above case, where $L=1$, $\varepsilon=0.1$.}
    \label{error}
    \enf
Meanwhile, Fig. \ref{error} shows the difference between $u$ and $u_\varepsilon$. The graph shows that away from the exit arc of small distance, $u_\varepsilon$ can approximate $u$ with small error of order $O(\varepsilon)$.

The second goal of the experiment is to see whether the Neumann-Robin model can be perfectly applied to solve the escape problem with different radius and different neck length.
\bef[ht]
\centering
    \includegraphics[width=2.5in]{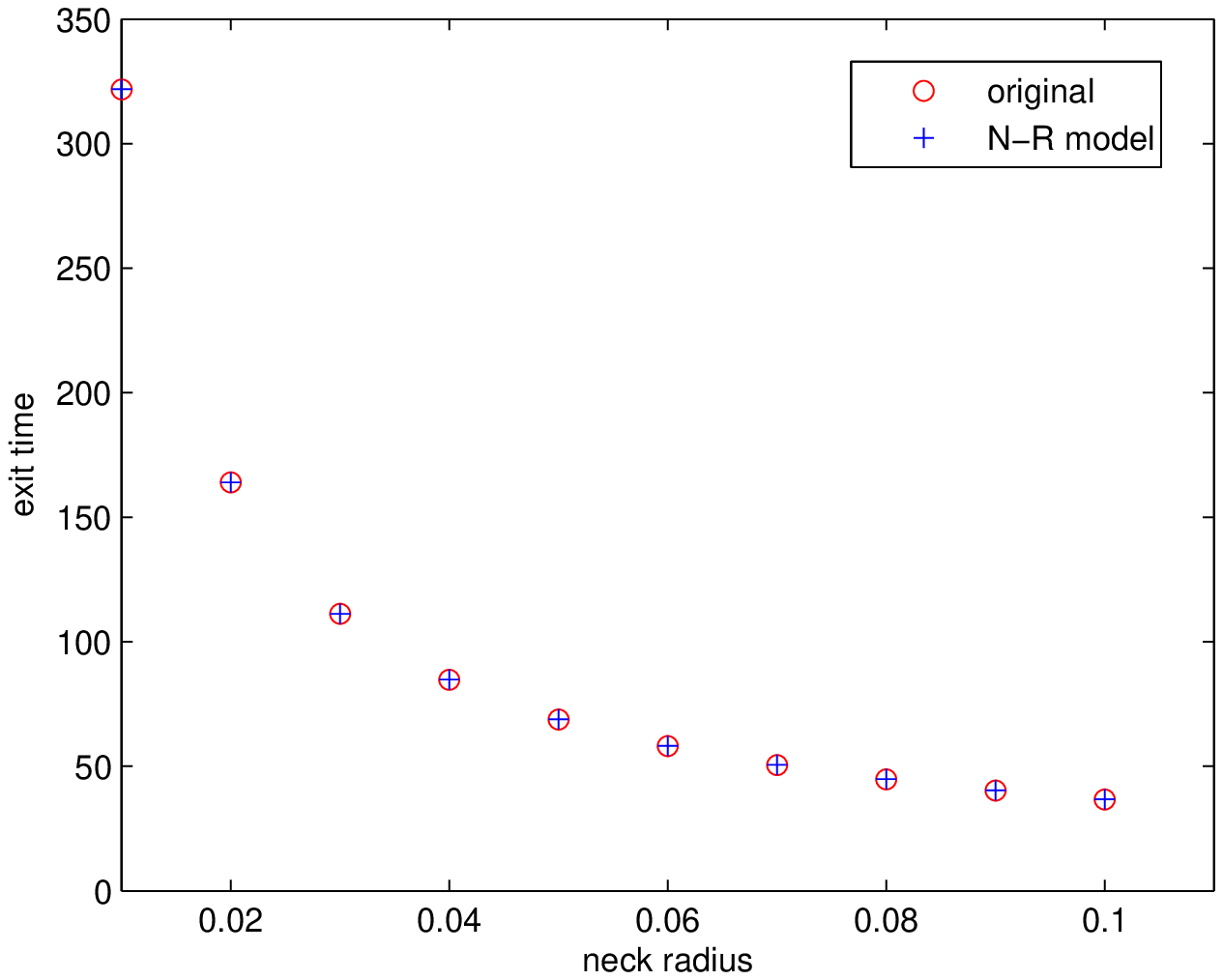} \includegraphics[width=2.5in]{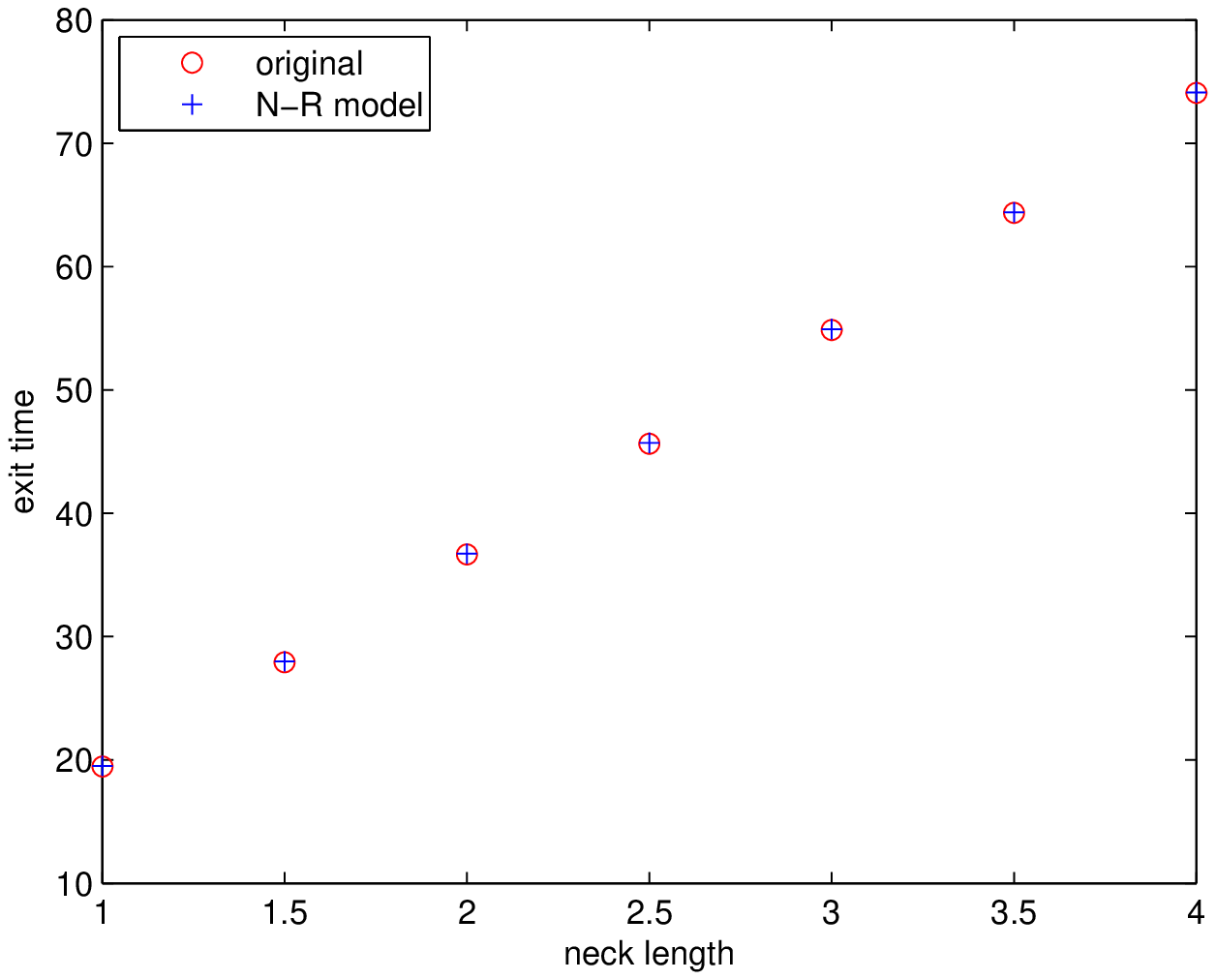}\\
    \caption{Left figure: Comparison between numerical results of original narrow escape problem (1.1) and expansion formula (4.6) derived from Neumann-Robin model, with different neck radius. We choose head to be unit disk, the neck length to be $L=2$, vary the neck radius $\varepsilon$ from 0.01 to 0.1. Right figure: Comparison between numerical results of original narrow escape problem (1.1) and asymptotic formula (4.6) derived from Neumann-Robin model, with different neck length. Fix head to be unit disk, radius of neck to be $\varepsilon=0.1$, vary neck length $L$ from 1 to 4.}
    \label{radius}
    \enf
From Fig. \ref{radius}, we can see the Neumann-Robin model perfectly solves the escape problem. The numerical solution of these two problems match with each other within error of order $O(\varepsilon)$. The figure on the left side is the case with different radius, 'o' represents the numerical solution of the original narrow escape problem, '+' represents the asymptotic formula (4.6) for Neumann-Robin model. Easy to see they are coincide. The figure on the right side is the case with different neck length. Similarly, the results coincide with each other.

\begin{table}
\centering
\begin{threeparttable}
\caption{Comparison result for center point of domain}
\label{tableu}
\begin{tabular}{llllllll}
\toprule
$\varepsilon$ & $L$ & $u_r$ & $u_\varepsilon$ & $u_\varepsilon-u_r$ & $u$ & $u-u_\varepsilon$ & $O(\varepsilon)$\\
\midrule
0.1	&1	  &19.4569	&19.5136   &-0.0515&19.5651	        &-0.0567	&0.1   \\
0.1	&1.5  &27.9232	&27.9914   &-0.0527&28.0441			&-0.0682	&0.1   \\
0.1 &2	  &36.6689	&36.7189   &-0.0542&36.7731           &-0.05  &0.1   \\
0.1	&2.5   &45.638	&45.6962    &-0.0558&45.752       	 &-0.0582   &0.1   \\
0.1	&3	  &54.8554  &54.9235    &-0.0575&54.981				&-0.0681&0.1   \\
0.1 &3.5	&64.354  &64.4007   &-0.0593&64.46	         &-0.0467	&0.1  \\
0.1	&4	  &74.0799  &74.1279     &-0.0611&74.189	     &-0.048  	&0.1    \\
0.09&2	  &40.2571  &40.3195    &-0.0493 &40.3688	     &-0.0624  &0.09  \\
0.08&2    &44.7399  &44.8107   &-0.0437 &44.8544	     &-0.0708   &0.08  \\
0.07&2	  &50.4918  &50.5589   &-0.0407 &50.5996	        &-0.0671&0.07  \\
0.06&2	  &58.1414  &58.1882   &-0.0323  &58.2205        &-0.0468	   &0.06 \\
0.05&2    &68.7906	&68.8463    &-0.0273 &68.8736	   	&-0.0557   &0.05  \\
0.04&2    &84.7263  &84.7901    &-0.0228   &84.8129		&-0.0638  &0.04  \\
0.03&2	  &111.2003  &111.2612   	&-0.0178&111.279	    &-0.0609   &0.03  \\
0.02&2    &	163.9653  &164.0201  &-0.0134&164.0335	     &-0.0548   &0.02  \\
0.01&2    &	321.7331  &321.8301 &-0.0134 &321.8435      &-0.0504	&0.01   \\
\bottomrule
\end{tabular}
\begin{tablenotes}
\small
\item $\varepsilon$: half length of the exit arc. $L$: length of spine neck. $u_r$: numerical solution of Neumann-Robin model. $u_\varepsilon$: value of asymptotic formula(4.6). $u_\varepsilon-u_r$: difference between $u_\varepsilon$ and $u_r$. $u$: numerical solution of narrow escape problem(1.1). $u_\varepsilon-u_r$: difference between $u_\varepsilon$ and $u_r$. $O(\varepsilon)$: error term.
\end{tablenotes}
\end{threeparttable}
\end{table}

Table.\ref{tableu} concretely demonstrates the comparison of these two problems. Releasing the particle at the center of the spine head, Table 5.1 shows the numerical results of the exact solution $u$, Neumann-Robin model solution $u_r$ and asymptotic formula $u_\varepsilon$, with respect to different neck radius $2\varepsilon$ and neck length $L$. First, from the comparison of the results in Neumann-Robin model $u_r$ and asymptotic formula $u_\varepsilon$, we can see our derivational calculation for Neumann-Robin problem by using layer potential techniques in section 3 is correct as shown in theorem 3.1. Second, from the comparison between $u$ and $u_\varepsilon$, we can see their difference is of order $O(\varepsilon)$, which means the asymptotic formula can more precisely approximate the exit time.
\\
\bef
\centering
  \includegraphics[width=1.8in]{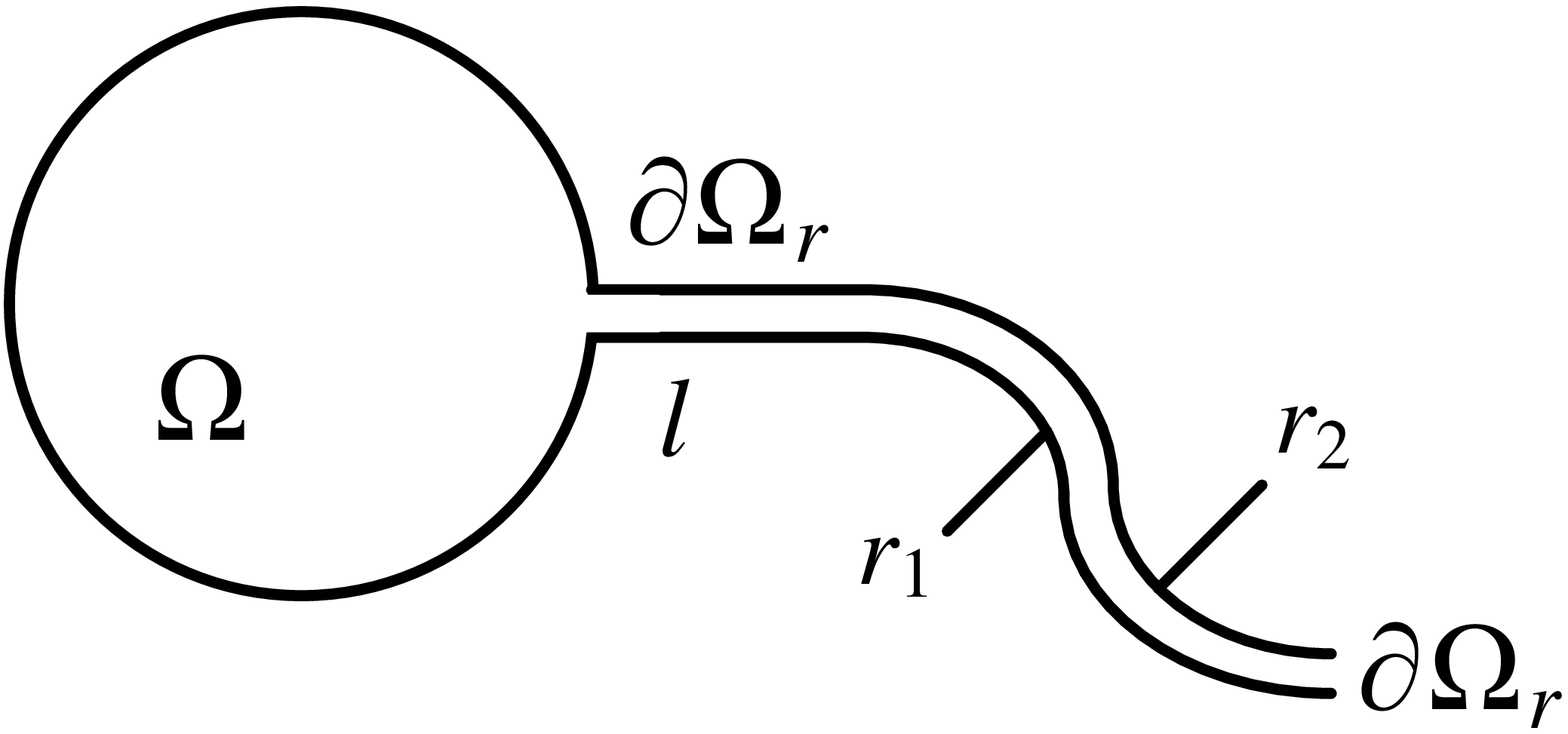}\\
  \caption{The approximated geometry model of dendritic spine with non-straight long spine neck, where $\Omega$ is the domain with long neck, $\partial \Omega_{r}$ is the reflection part, $\partial \Omega_{a}$ is the absorbing part. $l$ is the length of straight part of neck. The non-straight part of spine neck is composed of two circle arc, whose radius are $r_1$ and $r_2$ respectively.}
  \label{non}
 \enf

Next, the domain with non-straight spine neck is considered. Suppose that $\Omega$ is the spine domain with a smooth spine neck, $\Omega_{h}$ is the spine head, $\Gamma_{\varepsilon}$ is an arc of center $x^*$ which connects spine head and spine neck, $L$ is the absolute length of the spine neck, and $\kappa(x)$ is the curvature at the point $x$. The approximated geometry is the same as the previous situation, but with non-straight neck, see Fig. \ref{non}.

Through a number of numerical simulations, we find out that the Neumann-Robin model can be easily applied to solve the narrow escape problem in such a situation. But because of the curvature on the neck, we need to find a better corresponding neck length, not just the absolute value of the neck length $L$. We eventually see that if we insert
$$\tilde{L}=L+\varepsilon \int_{L}\kappa(x) ds$$
into formula (4.6), where  $L$ is the absolute length of the neck, $\kappa(x)$ is the curvature of the point $x$ on spine neck, then the Neumann-Robin model in section 3 can approximate the exit time even in the non-straight spine neck case. The first mean passage time $u_\varepsilon(x)$ of a Brownian particle confined in $\Omega_h$ exiting through $\partial\Omega_a$ in the non-straight spine neck domain (Fig.\ref{non}) can be approximated by the formula
\beq
u_\varepsilon(x)\approx\frac{|\Omega_{h}|\tilde{L}}{2\varepsilon}+\frac{|\Omega_{h}|}{\pi}(\frac{3}{2}+\ln\frac{1}{2\varepsilon})+ \frac{\tilde{L}^{2}}{2}+\Phi_{\Omega_{h}}(x,x^*),
\end{equation}
where $\Phi_{\Omega_{h}}(x,x^*)$ is the same as (3.24) for the domain $\Omega_h$. Our experimental data show that the error between this formula (5.1) and the exact solution for (1.1) in the non-straight spine neck domain is of order $O(\varepsilon)$.

There is one example. Consider the domain (Fig.\ref{non}) which is composed of unit disk head, the straight neck part  $l=1$, non-straight part $r_1=1$, $r_2=1$, and exit arc length $2\varepsilon$, $\varepsilon=0.1$. The numerical results of the solution $u$ for (1.1) and the result $u_\varepsilon$ of the expansion formula (5.1) with $\tilde{L}=L+\varepsilon \int_{L}\kappa(x) ds$ are given in Fig.\ref{co1111}.

\bef
\centering
  \includegraphics[width=6in]{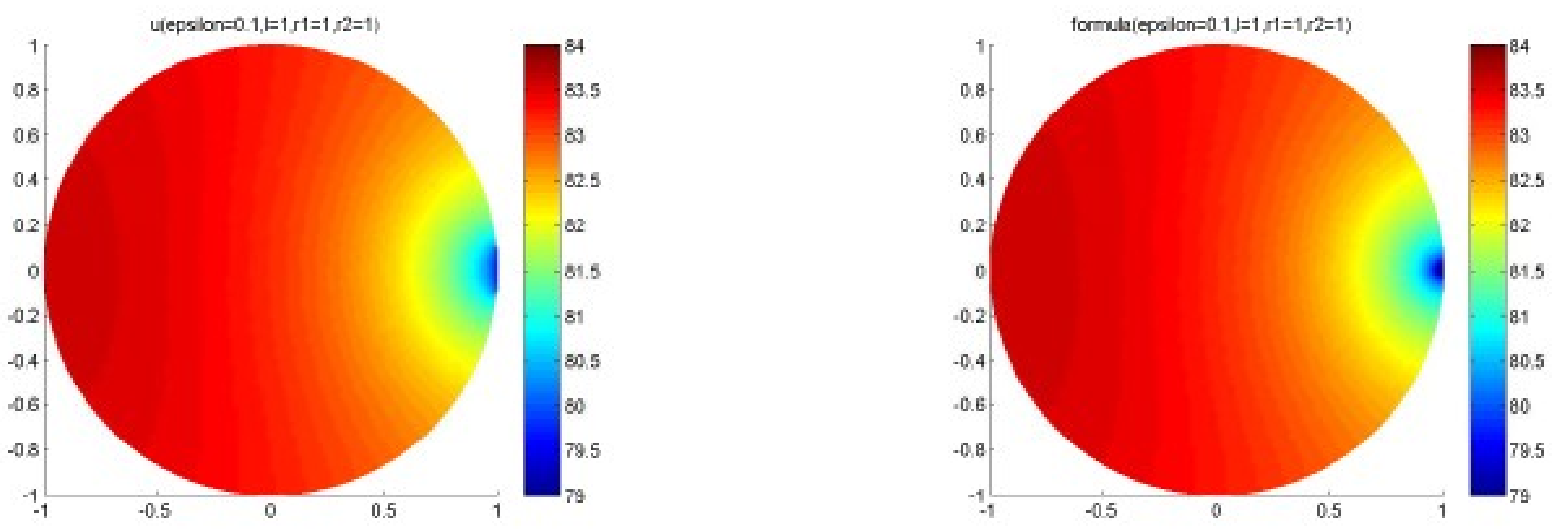}\\
  \caption{The spine geometry is approximated by a unit disk head, straight neck length $l=1$ and two quarters of unit disk with exit length $2\varepsilon=0.2$. Coordinates represent the position of the iron where it is initiated. The color represents the exit time of particle initiated at this point. Left figure: the numerical result $u$ of (1.1). Right figure: exit time computed by asymptotic formula (5.1).}
  \label{co1111}
 \enf
The figure on the left side gives the numerical solution $u$ of (1.1) in the spine head domain. The value at each point means the exit time of the particle initiated at that point. The figure on the right side is the solution $u_\varepsilon$ for asymptotic formula (5.1) with the neck length $\tilde{L}$. From these two figures, we can easily see the results in both situations agree within a small error. The numerical data show that the error between these two situations are of order $O(\varepsilon)$. Meanwhile, Fig. \ref{error2} shows the difference between $u$ and $u_\varepsilon$. The graph shows that away from the exit arc of small distance, $u_\varepsilon$ can approximate $u$ within small error of order $O(\varepsilon)$. The data of Table \ref{tableun} also confirm this assertion.

\bef
\centering
    \includegraphics[width=2.5in]{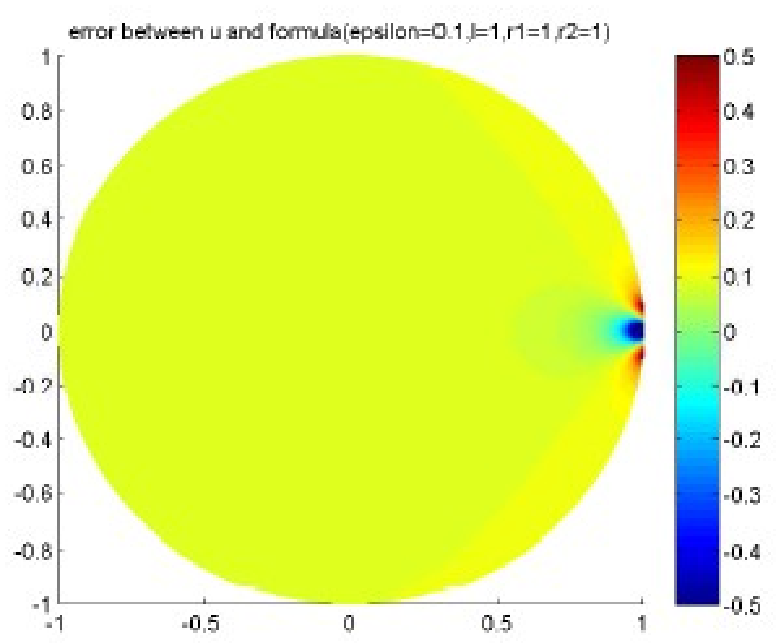} \\
    \caption{Error between $u$ and $u_\varepsilon$ in the above case, where $\varepsilon=0.1$, $l=1$, $r_1=1$, $r_2=1$.}
    \label{error2}
    \enf

\begin{table}
\centering
\begin{threeparttable}
\caption{Comparison result for center point of domain}
\label{tableun}
\begin{tabular}{llllllll}
\toprule
$\varepsilon$ & $l$ & $r_1$ & $r_2$ & $u$ & $u_\varepsilon$ &$u_\varepsilon-u$ & $O(\varepsilon)$\\
\midrule
0.1	&1  &0.7&0.9   &70.4851	&70.7957	&0.3106 &0.1\\
0.1	&1.5&0.7&0.9	&80.36	&80.6873	&0.3273&0.1\\
0.1	&2	&0.7&0.9	&90.5009	&90.8358	&0.3349&0.1\\
0.1	&1	&1  &1      &82.9631	&83.253	&0.2899&0.1\\
0.05&1	&1  &1	&148.1522	&148.3267	&0.1745&0.05\\
0.05&2	&1  &1	&184.3655	&184.5477	&0.1822&0.05\\
0.05&3	&1  &1	&221.5875	&221.7754	&0.1879&0.05\\
\bottomrule
\end{tabular}
\begin{tablenotes}
\small
\item $\varepsilon$: half length of the exit arc. $l$: length of straight part of spine neck. $r_1$ and $r_2$: the radius of two circles. $u_\varepsilon$: value of asymptotic formula(5.1). $u$: numerical solution of narrow escape problem(1.1). $u_\varepsilon-u$: difference between $u_\varepsilon$ and $u$. $O(\varepsilon)$: error term.
\end{tablenotes}
\end{threeparttable}
\end{table}

\section{Conclusion}
In this paper, using Neumann-Robin Boundary Model we transform spine singular domain to smooth spine head domain. We provided mathematically rigorous derivation of the leading order term in the asymptotic expansion of the solution of Neumann-Robin Boundary Model which we invented to solve narrow escape problem in domain with long neck. The result shows that using this model we find first escape time up to order $O(\varepsilon)$ which is not found in other papers. The solution to the Neumann-Robin Boundary Model in the spine head domain can be easily applied to the calcium diffusion model of the narrow escape problem, one with straight spine neck and the other with non-straight spine neck. As for the non-straight spine neck, integrating the neck curvature we can get an effective length which can be put into our explicit expansion formula of Neumann-Robin Boundary Model and get the approximated exit time. This Neumann-Robin Model can be extended to three dimension. Because of the existence of long neck the operators we defined in section 3 can be shown bounded and expansion formula can be similarly derived with two dimensional case. What's more, for narrow escape problem in three dimensional smooth domain with exit on the boundary, since the difficulty in 3D case, up to now, only spherical cases are appeared in other papers. But for any other smooth domain in 3D, even it has no long neck, we can still use this Neumann-Robin Boundary Model by punching long neck to the small exit on the boundary which will be the subject of a forthcoming paper.

\section{Acknowledgement}
This work was funded by Korean National Research Foundation through NRF gtants Nos. 2009-0085987 and BK21+ at Inha University. Thanks for Hyundae lee and Hyeonbae Kang's fruitful discussion and suggestion.

\end{document}